\documentclass{amsart}

\usepackage{amssymb}
\usepackage[all]{xy}
\usepackage{hyperref}

\usepackage{enumitem}   

\setlist[enumerate]{itemsep=.2em,topsep=.2em,leftmargin=1.25em,itemindent=2.0em}

\newtheorem{thm}{Theorem}

\newtheorem{cor}[thm]{Corollary}

\newtheorem{prop}[thm]{Proposition}

\newtheorem{conj}[thm]{Conjecture}

\theoremstyle{definition}

\newtheorem{say}[thm]{}

\newtheorem{rem}[thm]{Remark}
\newtheorem{rems}[thm]{Remarks} 

\newtheorem*{ack}{Acknowledgments}      

\newtheorem{defn-thm}[thm]{Definition--Theorem}  
\newtheorem{defn-lem}[thm]{Definition--Lemma}  

\theoremstyle{remark}

\renewcommand{\c}[0]{{\mathbb C}}  

\renewcommand{\o}[0]{{\mathcal O}} 
\newcommand{\z}[0]{{\mathbb Z}}

\newcommand{\p}[0]{{\mathbb P}}

\newcommand{\map}[0]{\dasharrow}
\newcommand{\qtq}[1]{\quad\mbox{#1}\quad}

\newcommand{\pic}[0]{\operatorname{Pic}}

\newcommand{\pico}[0]{\operatorname{\mathbf{Pic}}^{\circ}}

\newcommand{\ns}[0]{\operatorname{NS}}

\newcommand{\supp}[0]{\operatorname{Supp}}

\newcommand{\jac}[0]{\operatorname{\mathbf{Jac}}}

\newcommand{\grass}[0]{\operatorname{Grass}}

\newcommand{\tsum}[0]{\textstyle{\sum}}

\def\into{\DOTSB\lhook\joinrel\to}

\def\loccoh#1.#2.#3.#4.{H^{#1}_{#2}(#3,#4)}

\DeclareMathAlphabet{\mathchanc}{OT1}{pzc}%
                                {m}{it}

\usepackage[all]{xy}\xyoption{dvips}

\newcommand{\mw}[0]{\operatorname{MW}}

\begin{document}
\bibliographystyle{amsalpha}

 \hfill\today

 \title[Sections of Jacobian fibrations]{Sections of Jacobian fibrations over lines}
 
 \author{J\'anos Koll\'ar and Giulia Sacc\`a}

 \begin{abstract}  Let $|H|$ be a linear system on a smooth surface $S$.
   We study the cohomology classes  of  sections of the universal 
   Jacobian  over lines in $|H|$. When $S$ is a K3 surface, 
the universal  compactified
   Jacobian is  a hyperk\"ahler manifold, and 
we do not see how to reconcile our results  with some of the steps in  [BKV25].
       \end{abstract}

 \maketitle

 Let $p:X\to \p^n$ be a Lagrangian fibration of a  smooth, hyperk\"ahler variety, and 
  $L\subset \p^n$  a general line. Assume that there is a section
 $\sigma_L: L\to X_L:=p^{-1}(L)$.  \cite{bkv2025} studies whether
 it can be extended to a meromorphic section $\sigma_P:\p^n\map X$,
 possibly after changing  the complex structure of $X$.

 We show that in many cases the answer is negative; there are topological obstructions. The right setting seems to be universal Jacobians for linear systems of curves. We work over $\c$.

 \begin{prop}\label{main.S.exmp}
 Let $S$ be a  smooth, projective surface 
 such that $\pic(S)=\z[H]$, where $|H|$  is basepoint-free,  and   members of $|H|$ have at worst nodes in codimension 1 on $|H|$. 
  Let $p: J(S, H)\to  |H|$ be the universal compactified Jacobian as in 
 (\ref{univ.J.say}.2),   
   $L\subset |H|$ 
   a general line,  $J_L:=p^{-1}(L)$, and $g$  the genus of the curves in $|H|$. Then
   \begin{enumerate}
   \item  the Mordell-Weil group of $J_L\to L$ is isomorphic to $\z^{r-1}$ for $r=(H^2)$.
     \item Let   $Z\subset J_L$ be a section whose cohomology class is
    contained in the image of  the restriction map
    $H^{2g}(J(S, H), \z)\to H^{2g}(J_L, \z)$. Then $Z$ is the zero section.
    \end{enumerate}
 \end{prop}

\begin{cor} \label{main.S.exmp.cor} Using the notation of Proposition~\ref{main.S.exmp},
  let $\sigma_L:L\to J_L$ be a holomorphic section. Assume that there is some complex structure $\mathchanc{J}$ on $X:=J(S, H)$ such that
  $p:(X,\mathchanc{J})\to |H|$ is holomorphic, and $\sigma_L$ 
  extends to a meromorphic section $\sigma:|H|\map (X, \mathchanc{J})$.
  Then $\sigma_L$ is the zero section.
\end{cor}

Proof. Let $Z\subset (X,\mathchanc{J})$ be the closure of the image of
$\sigma$. Since $X$ is smooth, it defines a cohomology class $\eta_P\in H^{2g}(X, \z)$ whose restriction to $J_L$ is  $\eta_L\in H^{2g}(J_L, \z)$, the
cohomology class of $Z_L:=\sigma_L(L)$. \qed

\medskip

\begin{rems} If $S$ is a K3 surface, 
   Corollary~\ref{main.S.exmp.cor} seems inconsistent with some of the intermediate claims in \cite{bkv2025}, but not with the main theorem.
  
 Other results creating sections using Tate-Shafarevich twists
  are discussed in \cite{k-abtwist}.
  Twists of Lagrangian fibrations are also treated in
  \cite{MR4925350, abasheva2024shafarevichtategroupsholomorphiclagrangian, sacca2025}.

 See Proposition~\ref{main.S.exmp.2} for another version of  Proposition~\ref{main.S.exmp}.
  \end{rems}

Before we start the proof of 
Proposition~\ref{main.S.exmp}, let us recall some facts about families of Jacobians.

 \begin{say}[Universal Jacobian] \label{univ.J.say}
   Let $S$ be a smooth, projective  surface, and  $|H|$ a basepoint-free linear system on  $S$.
   We set $d:=\dim |H|$, $r:=(H^2)$ and let $g$ denote the genus of the curves in $H$.
   
Let $|H|^\circ\subset |H|$ be the open subset parametrizing irreducible curves with at worst nodes. 
The   {\it universal compactified Jacobian}
$$
p^\circ: J^\circ=J^\circ(S, H)\to |H|^\circ
 \eqno{(\ref{univ.J.say}.1)}
$$
 has the following properties; see \cite{MR569548} or \cite[Sec.B]{MR1658220} for details. 
 
For a smooth  curve $C\in |H|$, the fiber of $p^\circ$ over
$[C]$ is  $\jac(C)$. If $C$ is irreducible but singular,
then the fiber is the compactification of $\jac(C)$.
Then $p^\circ$ is flat, projective, and 
$J^\circ(S, H) $ is smooth.

The trivial line bundle on each $C$ provides a section  $|H|^\circ\to J^\circ$.

If $\pic(S)=\z[H]$, then   $p^\circ$ extends to a
flat, projective morphism
$$
p: J=J(S, H)\to |H|,
 \eqno{(\ref{univ.J.say}.2)}
 $$
 but
 $ J(S, H)$ is not smooth in general.
 However $ J(S, H)$ is smooth if $S$ is a K3 surface and $\pic(S)=\z[H]$ by 
 \cite {MR751133}.
 \end{say}

 \begin{say}[Sections over lines] \label{K3.LF.exmp}
Assume in addition that   $|H|\setminus |H|^\circ$ has codimension $\geq 2$.

 Let
 $L\subset |H|^\circ$ be a   line and  $|H_L|$ the corresponding   pencil.
 Assume that it has $r:=(H^2)$ distinct basepoints
  $\{c_i:i=1,\dots, r\}$.  Blowing them up we get  a basepoint-free pencil
  $B_{r}S\to L $. The $r$  exceptional curves give $r$ sections
 $E_i$  of $B_{r}S\to L$.

   Let $M$ be a line bundle on $B_rS$ such that $(L\cdot F)=0$, where $F$ is a general fiber. Restricting $M$ to each fiber of  defines a section  $Z_M\subset J_L$.
   
   The {\it Shioda-Tate isomorphism} says that
   if $h^1(S, \o_S)=0$  and all members of  $|H_L|$ are irreducible, then $M\mapsto Z_M$
    gives an isomorphism between the
   orthogonal complement of $F$ in the N\'eron-Severi group $\ns(B_rS)/\langle F \rangle$,
   and the  Mordell-Weil group $\mw(J_L/L)$ of $J_L\to L$,  see
   \cite{MR1202625, MR1610977, MR1714831, MR2494474}.

   We will focus on sections constructed from the
   exceptional curves $E_i$. 
 For $m_i\in \z$ sstisfying $\tsum m_i=0$, the divisor
   $\tsum m_iE_i$   gives a section of $J_L\to L$.
 For $\mathbf{m}:=(m_1, \dots, m_{r})$ let us denote 
   the  corresponding section of $J_L\to L$  by
   $$
   \begin{array}{ll}
     \sigma_L(\mathbf{m})&\mbox{as a morphism $L\to J_L$,}\\
    Z_L(\mathbf{m})&\mbox{as a subvariety of $ J_L$, and}\\
    \eta_L(\mathbf{m})&\mbox{as a class in $ H^{2g}(J_L, \z)$.}
    \end{array}
   \eqno{(\ref{K3.LF.exmp}.1)}
   $$
   These give a subgroup  $\Gamma_L\subset \mw(J_L/L)$.
   Note that $\Gamma_L\cong \z^{r-1}$, and
    $\Gamma_L=\mw(J_L/L)$ iff $\pic(S)=\z[H]$. 
   
   Let $W\subset \grass(\p^1, |H|)$ be the open subset parametrizing such pencils.
    Over $W$, the 
    $(m_1, \dots, m_{r})$ are determined only up to then action of  the monodromy group.  
   In the universal family over $W$,   the monodromy group acts transitively on the basepoints.  (In most cases,  it is  the full symmetric group.)
 \end{say}

 \begin{say}[Proof of Proposition~\ref{main.S.exmp}]
   \label{main.S.exmp.pf}
   As we noted in Paragraph~\ref{K3.LF.exmp},
   the Shioda-Tate isomorphism implies (\ref{main.S.exmp}.1).
   This was already noted in \cite{MR1714831}.

   As in (\ref{K3.LF.exmp}.1), let
    $\eta_L(m_1, \dots, m_{r})\in H^{2g}(J_L, \z)$ be a  cohomology class.
   As we move the line $L\subset |H|$ in the Grassmannian of lines,
   the $m_i$    get permuted.
   By contrast, the elements of the image of
   $H^{2g}(J, \z)\to H^{2g}(J_L, \z)$
   are invariant under monodromy.
   Assume now that
   $$
   \eta_L(m_1, \dots, m_{r})=\eta_L(m'_1, \dots, m'_{r})
   \qtq{iff} m_i=m'_i \quad \forall i.
   \eqno{(\ref{main.S.exmp.pf}.1)}
   $$
Then $\eta_L(m_1, \dots, m_{r})$ is monodromy invariant  only  when
   $m_1= \cdots = m_{r}$, but in $\mw(J_L/L)$ we have $\sum m_i=0$.
   This shows (\ref{main.S.exmp}.2), provided (\ref{main.S.exmp.pf}.1) holds.\qed
 \end{say}

   By the Shioda-Tate isomorphism the sections
   $Z_L(m_1, \dots, m_{r})$ are all different from each other.
   Thus we are done if the map from sections to homology classes is an injection.  This is what we discuss next.

   We start with the strongest possible form, but then we prove only a rather special case, which however implies (\ref{main.S.exmp.pf}.1).

 \begin{conj}\label{mw.t.h2.conj} Let $C$ be a smooth, projective curve over a  field $k$, 
   and $p:A\to C$ a smooth group scheme with irreducible fibers, whose generic fiber $A_K$ is an Abelian variety. Let $Z_0$ be the zero section, and assume that $A_K$ has no Abelian subvarieties defined over $k$. 
  Then  sending a section to its numerical equivalence class
   $$
    \mw(A/C)\ni Z\mapsto [Z]-[Z_0]\in  N_1(A)
   \eqno{(\ref{mw.t.h2.conj}.1)}
   $$
   is an injective   group homomorphism modulo torsion.
   \end{conj}

 \begin{say}[Special case of Conjecture~\ref{mw.t.h2.conj}]\label{mw.to.h2.say}
   Let $T$ be a smooth, projective surface over a field $k$ such that $h^1(T, \o_T)=0$, and  $q:T\to \p^1$ a morphism  whose fibers are irreducible and nodal.
   Let $J_T\to \p^1$ be the relative Jacobian.
   Fix a section  $E_0$  of $q$ and let $F$ be a general fiber of $q$.

    Let $K$ be the function field of $\p^1$ and  $(T_K, e_K)$ the generic fiber with the point given by $E_0$.
 This defines  embeddings
 $$
 j_K:T_K\into \jac(T_K)\qtq{and}  j: T\into J_T.
 \eqno{(\ref{mw.to.h2.say}.1)}
 $$
 Note that pull-back induces an isomorphism
 $$
 j_K^*: \pico\bigl(\jac(T_K)\bigr)\cong \pico(T_K)=\jac(T_K).
 \eqno{(\ref{mw.to.h2.say}.2)}
 $$
 Therefore, given a  line bundle $M$ on $T$ such that $(M\cdot F)=0$,
 there is a line bundle $M_J$ on $J_T$ such that
 $$
 j^*M_J= M \qtq{as elements of}  \ns(T)/\langle F \rangle.
 \eqno{(\ref{mw.to.h2.say}.3)}
 $$
 Let $D$ be a divisor on $T$ such that $(D\cdot F)_T=0$.  Then
  $$
 (D\cdot M)_T=(D\cdot j^*M_J)_T=(j_*(D)\cdot M_J)_{J_T}.
 \eqno{(\ref{mw.to.h2.say}.4)}
 $$
 The middle number is  well defined since $(D\cdot F)_T=0$.

 (Note that $(D\cdot M)_T$  is the product  used in the construction of the   Shioda height pairing  \cite{MR1202625, MR1714831}. We need to lift it to the Jacobian.)

 Let now $E_0,\dots, E_r$ be sections of $q$.
 Then $D_i:=E_i-E_0$ are linearly independent elements of the
 orthogonal complement of $F$ in the N\'eron-Severi group $\ns(T)/\langle F \rangle$.
 By   the Hodge index theorem, 
 the intersection form  is  negative definite on
 $\langle D_1, \dots, D_r \rangle$.

 Combining with (\ref{mw.to.h2.say}.4) we conclude that
 $$
 j_*: \langle D_1, \dots, D_r \rangle\to H_2(J_T, \z)
 \qtq{is injective.}
 \eqno{(\ref{mw.to.h2.say}.5)}
 $$
 This is almost what we want, except that the elements of
 $j_*\langle D_1, \dots, D_r \rangle$ are linear combinations of sections
 of $T\to \p^1$, viewed as section of $J_T\to \p^1$.
 For (\ref{main.S.exmp.pf}.1) we need their sum in the Mordell-Weil group.
 That is, if $Z_1, Z_2, Z_3$ are sections
 such that $Z_3=Z_1\oplus Z_2$ in the Mordell-Weil group, then
 $$
 (Z_1\cdot M_J)+ (Z_2\cdot M_J)=(Z_3\cdot M_J)+(Z_0\cdot M_J).
 \eqno{(\ref{mw.to.h2.say}.6)}
 $$
 Let $\tau_i:A\to A$ denote translation by $Z_i$.  Then
 $(Z_i\cdot M_J)=(Z_0\cdot \tau_i^*M_J)$, hence
 (\ref{mw.to.h2.say}.6) becomes
 $$
 (Z_0\cdot \tau_1^*M_J)+ (Z_0\cdot\tau_2^* M_J)=(Z_0\cdot \tau_3^*M_J)+(Z_0\cdot M_J).
 \eqno{(\ref{mw.to.h2.say}.7)}
 $$
 This holds since $\tau_1^*M_J\otimes \tau_2^*M_J\cong M_J\otimes \tau_3^*M_J$ by the Theorem of the square
 \cite[p.60]{mumf-abvar}.
 \end{say}

 The proof also shows the following variant of
 Proposition~\ref{main.S.exmp}.

 \begin{prop}\label{main.S.exmp.2}
   Let $S$ be a  smooth, projective surface  such that $h^1(S, \o_S)=0$,
   and 
 $|H|$  a basepoint-free linear system whose   members  are irreducible and  nodal in codimension 1 on $|H|$.
 
   Then the Mordell-Weil group of $p^\circ: J^\circ(S, H)\to  |H|^\circ$
   is isomorphic to $[H]^\perp$, the orthogonal complement of
   $[H]$ in $\ns(S)$.  \qed
 \end{prop}

 \begin{rem}  In particular,  every rational section of
   $p: J(S, H)\to  |H|$ is regular.
   We do not know whether the same holds for all
   Lagrangian fibrations of   smooth, hyperk\"ahler varieties.
 \end{rem}

 \begin{say}[Continuous sections]
   Let $\eta\in H_2(S, \z)$ be a homology class
   such that $c_1(H)\cap \eta=0$.
  Let  $N:=\sum n_iN_i$ be a representative of $\eta$, where the $N_i\subset S$ are oriented submanifolds. Under suitable transversality assumptions,
   $ C\mapsto  \o_C(N\cap C)$ defines  a continuous section $\sigma_N$
   of $p:J(S, H)\to |H|$.
   If the $N_i$ are algebraic, then $\sigma_N$ is an algebraic section.
   (There are some problems when a singular point of $C\in |H|$ is in  $\supp N$.)

   If $S$ is a K3 surface, then every homology class is algebraic
   on some deformation of $S$, so these continuous sections become homotopic to algebraic sections, after changing the complex structure.

   Continuous sections of  $J_L\to L$ come from
   homology classes of
   $$
  H_2(B_rS, \z)\cong  H_2(S, \z)\oplus \langle [E_1],\dots, [E_r]\rangle.
   $$

Since our obstructions are topological, they apply to all
deformations of $J\to \p^d$, not just to those  that come from a deformation of $(S, H)$.
Such an example  is given by the O'Grady 10-folds $X\to \p^5$ \cite{MR1703077}, which deform to 
 LSV 10-folds \cite{MR3710794}. 

   \end{say}

 \begin{ack}  We thank  F.~Bogomolov, L.~Kamenova, and M.~Verbitsky  for  a series of  e-mails. These answered several of our questions, but we failed to converge to an agreement on all the points.

   Partial  financial support  to JK   was provided  by the Simons Foundation   (grant number SFI-MPS-MOV-00006719-02).
Partial  financial support  to GS   was provided  by the NSF
    (grant number DMS-2144483) and the Simons Foundation   (grant number SFI-MPS-MOV-00006719-09).

   \end{ack}


 \def\cprime{$'$} \def\cprime{$'$} \def\cprime{$'$} \def\cprime{$'$}
  \def\cprime{$'$} \def\dbar{\leavevmode\hbox to 0pt{\hskip.2ex
  \accent"16\hss}d} \def\cprime{$'$} \def\cprime{$'$}
  \def\polhk#1{\setbox0=\hbox{#1}{\ooalign{\hidewidth
  \lower1.5ex\hbox{`}\hidewidth\crcr\unhbox0}}} \def\cprime{$'$}
  \def\cprime{$'$} \def\cprime{$'$} \def\cprime{$'$}
  \def\polhk#1{\setbox0=\hbox{#1}{\ooalign{\hidewidth
  \lower1.5ex\hbox{`}\hidewidth\crcr\unhbox0}}} \def\cdprime{$''$}
  \def\cprime{$'$} \def\cprime{$'$} \def\cprime{$'$} \def\cprime{$'$}
\providecommand{\bysame}{\leavevmode\hbox to3em{\hrulefill}\thinspace}
\providecommand{\MR}{\relax\ifhmode\unskip\space\fi MR }
\providecommand{\MRhref}[2]{%
  \href{http://www.ams.org/mathscinet-getitem?mr=#1}{#2}
}
\providecommand{\href}[2]{#2}

 Department of Mathematics, Princeton University,

 Fine   Hall, Washington Road, Princeton, NJ 08544-1000, USA

   \email{kollar@math.princeton.edu}

   \medskip

   Department of Mathematics, Columbia University 

   2990 Broadway, New York,  NY 10027, USA
   
\email{gs3032@columbia.edu}

\end{document}